\documentclass{amsart}
\usepackage[all]{xy}
\usepackage{amssymb,amscd}  %% can also include 'bm' for 'boldmath', used
\newtheorem{theorem}{Theorem}
\newtheorem{proposition}[theorem]{Proposition}
\newtheorem{lemma}[theorem]{Lemma}

\theoremstyle{definition}

\newtheorem{example}[theorem]{Example}
\newtheorem{remark}[theorem]{Remark} %% see below
\newcommand\Z{\mathbf{Z}}
\newcommand\tensor{\otimes}
\newcommand\isomorphic{\cong}

\DeclareMathOperator{\Hom}{Hom}

\newcommand\Directsum{\bigoplus}
\newcommand\directsum{\oplus}

\newcommand\union{\cup}

\newcommand\intersect{\cap}

\newcommand\abs[1]{{\left|#1\right|}}

\newcommand\w[1]{\wedge^{#1}V}
\newcommand\s[1]{S^{#1}V}
\newcommand\F{\mathrm{Fil}}
\newcommand\G{\mathrm{Gr}}
\newcommand\eps{\varepsilon}
\newcommand\dee{\partial}
\newcommand\Lbar{\overline{\Lambda}}
\newcommand\lbar{\overline{\lambda}}
\begin{document}

\title{On inverting the Koszul complex}

\author{Kamal Khuri-Makdisi}
\address{Mathematics Department and Center for Advanced Mathematical Sciences,
American University of Beirut, Bliss Street, Beirut, Lebanon}
\email{kmakdisi@aub.edu.lb}
\thanks{April 26, 2007}

\begin{abstract}
Let $V$ be an $n$-dimensional vector space.  We give a direct
construction of an exact sequence that gives a $GL(V)$-equivariant
``resolution'' of each symmetric power $S^t V$ in terms of
direct sums of tensor products of the form 
$\wedge^{i_1} V \otimes \dots \otimes \wedge^{i_p} V$.  This exact sequence 
corresponds to inverting the relation in the representation ring of $GL(V)$
that is described by the Koszul complex, and has appeared before in
work by B. Totaro, analogously to a construction of K. Akin involving
the normalized bar resolution.  Our approach yields a concrete
description of the differentials, and provides an alternate
direct proof that $\mathrm{Ext}^t_{\wedge (V^*)}(k,k) = S^t(V)$.
\end{abstract}

\maketitle

%\section{Introduction}
%\label{section1}

Let $k$ be a field, and let $V$ be an $n$-dimensional vector space over $k$
with $n \geq 2$.  We view the alternating and symmetric powers $\w{i}$ and
$\s{t}$ of $V$ as representations of $G = GL(V) \isomorphic GL(n, k)$; we
allow any $i,t \in \Z$, with the understanding that $\w{i}$ and $\s{t}$ are
zero unless $0 \leq i \leq n$ or $t \geq 0$.  Working in the
representation ring of $G$ (i.e., the Grothendieck group of algebraic
representations of $G$), we can write $\s{t}$ as a polynomial in the
fundamental representations $\w{1} = V, \w{2}, \dots, \w{n}$; as we shall
see below, the terms of this polynomial can be ordered in a natural way
with alternating signs.  
This suggests that there should exist an exact sequence of
representations of $G$ that concretely realizes this polynomial
expression for $\s{t}$.  This sequence, given in
equation~\eqref{equation2} below, has appeared in Sections~2
and~4 of~\cite{Totaro}, as a modification of a construction
of~\cite{Akin}, via the normalized bar resolution of $k$ as a module
for the exterior algebra $\wedge (V^*)$.  It is known that there is a natural
isomorphism 
\begin{equation}
\label{equation0}
  \mathrm{Ext}^t_{\wedge(V^*)}(k,k) = \s{t},
\end{equation}
and computing the Ext group using the normalized bar resolution as
in~\cite{Totaro} produces the exact sequence~\eqref{equation2}.  It is
however not immediate to write down the differentials explicitly, as
this involves chasing through various dualizations and natural
isomorphisms.  Moreover, one needs to have prior knowledge
of~\eqref{equation0} to use this approach.

In this note, we give a direct, self-contained construction of the
exact sequence~\eqref{equation2}, with explicit differentials, and
without using~\eqref{equation0}.  We give a straightforward proof that
the sequence is exact; this involves a total induction on $t$ that
uses the Koszul complex to analyze the sequence~\eqref{equation2} in
terms of the previous cases.  Our proof generalizes the argument
of~\cite{KKMsympowers} from $SL(2)$ to $GL(n)$.  We hope that the
structure of our argument, as well as our identification
of~\eqref{equation2} as an inversion of the Koszul complex, will be of
interest in other contexts.  Combining our proof of exactness with a
separate calculation to show that our sequence is the same as the one
coming from the normalized bar resolution, we obtain an independent
proof of~\eqref{equation0}. 

The first few cases of the exact sequence, for $1 \leq t \leq 4$, are
\begin{equation}
\label{equation1}
\begin{split}
& 0 \to V \to \s{1} \to 0,
  \qquad\qquad 
  0 \to \w{2} \to V\tensor V \to \s{2} \to 0,\\
& 0 \to \w{3}
    \to \bigl(V \tensor \w{2} \bigr) \directsum \bigl(\w{2} \tensor V\bigr)
    \to V \tensor V \tensor V \to \s{3} \to 0,\\
& 0 \to \w{4}
    \to \bigl(V \tensor \w{3} \bigr)
          \directsum \bigl( \w{2} \tensor \w{2} \bigr)
          \directsum \bigl( \w{3} \tensor V \bigr) \to \\
& \qquad
    \to \bigl(V \tensor V \tensor \w{2} \bigr)
          \directsum \bigl( V \tensor \w{2} \tensor V \bigr)
          \directsum \bigl( \w{2} \tensor V \tensor V \bigr) \to \\
& \qquad
    \to V \tensor V \tensor V \tensor V
    \to \s{4} \to 0. \\
\end{split}
\end{equation}  

In general, the sequence has the form
\begin{equation}
\label{equation2}
0 \to T_t^1 \xrightarrow{\delta} T_t^2 \xrightarrow{\delta} \cdots
  \xrightarrow{\delta} T_t^t \xrightarrow{\pi} \s{t} \to 0,
\end{equation}
for suitable differentials $\delta$ and $\pi$, where, for
$1 \leq p \leq t$, the $p$th term $T_t^p$ is given by 
\begin{equation}
\label{equation3}
\begin{split}
T_t^p &= \Directsum_{\substack{1 \leq i_1, \dots, i_p \leq n \\
                              i_1 + \dots + i_p = t\\}}
             \bigl( \w{i_1} \tensor \cdots \tensor \w{i_p} \bigr) \\
 & \isomorphic
    \Directsum_{\substack{\ell_1, \dots, \ell_n \geq 0 \\
                          \ell_1 + 2\ell_2 + \dots + n \ell_n = t \\
                          \ell_1 + \ell_2 + \dots + \ell_n = p\\}}
       \frac{(\ell_1 + \dots + \ell_n)!}{\ell_1! \cdots \ell_n!}
     \bigl[
        V^{\tensor \ell_1} \tensor 
           (\w{2})^{\tensor \ell_2}\cdots \tensor (\w{n})^{\tensor \ell_n}
     \bigr].\\
\end{split}
\end{equation}
We shall see in the proof of Theorem~\ref{theorem2} that it is reasonable
to allow $p=0$, provided we define $T_0^0 = k$ and $T_t^0 = 0$ for $t \neq
0$.  Note also that $T_t^t = V^{\tensor t}$ and that $T_t^1 = \w{t}$, which
may be zero (our conventions imply that $T_t^p = 0$ unless 
$p \leq t \leq np$; note also that the condition $i_1, \dots, i_p \leq n$
in \eqref{equation3} is redundant since otherwise $\w{i} = 0$.)
Our choice of letters $t$ and $p$ refers to the ``total degree'' 
(i.e., the effect of a scalar matrix in $G$) and ``partial degree'' (i.e.,
the number of parts) of a decomposable tensor
$\alpha_1 \tensor \cdots \tensor \alpha_p \in 
\w{i_1} \tensor \cdots \tensor \w{i_p} \subset T_t^p$.
The multinomial coefficient
$\frac{(\ell_1 + \dots + \ell_n)!}{\ell_1! \cdots \ell_n!}$
refers to a direct sum of several copies of the representation
$V^{\tensor \ell_1} \tensor \cdots \tensor (\w{n})^{\tensor \ell_n}$
of $G$.  
In the special case $n=2$, we see that $T_t^{t-i} \isomorphic 
\binom{t-i}{i} \bigl[V^{\tensor (t-2i)} \tensor (\w{2})^{\tensor i}\bigr]$,
and we recover the result of~\cite{KKMsympowers}.

Before describing the differentials in~\eqref{equation2}, we pause to
explain why 
\begin{equation}
\label{equation4}
T_t^1 - T_t^2 + \dots + (-1)^{t-1} T_t^t + (-1)^t \s{t} = 0
\end{equation}
in the representation ring of $G$.  This gives the polynomial expression of
$\s{t}$ in terms of the $\w{i}$, and arises from inverting the relation
between symmetric and alternating powers that is expressed by the exactness
of the Koszul complex
\begin{equation}
\label{equation5}
0 \to \w{n} \tensor \s{*} \to \dots
  \to \w{2} \tensor \s{*} \to V \tensor \s{*} \to \s{*} \to k \to 0.
\end{equation}
Here $\s{*} = \directsum_{t \geq 0} \s{t}$
is the symmetric algebra on $V$; it is naturally isomorphic to
the polynomial algebra $k[e_1, \dots, e_n]$, with $\{e_1, \dots, e_n\}$ a
basis for $V$.  The last term $k$ is the trivial representation of $G$,
i.e., the unit element of the representation ring of $G$; this term should
be viewed as $k[e_1, \dots, e_n]/\langle e_1, \dots, e_n\rangle$.

The morphisms in~\eqref{equation5}
(see for instance Section~XXI.4 of~\cite{LangAlgebra}
or Section~VII.2 of~\cite{MacLane}) 
increase the degree in each $\s{*}$  
component by one, so by taking $G$-equivariant Hilbert series with a formal
parameter $x$, we obtain the identity of formal power series in the
representation ring of $G$:
\begin{equation}
\label{equation6}
  \bigl(\sum_{t=0}^\infty \s{t}\cdot x^t
  \bigr)
  \bigl(1 - V \cdot x + \w{2} \cdot x^2
         - \dots + (-1)^n \w{n} \cdot x^n
  \bigr)
= 1.
\end{equation}
Now replace $x$ by $-x$ and invert the sum over the $\w{i}$ to yield
\begin{equation}
\label{equation7}
\begin{split}
\sum_{t=0}^\infty (-1)^t \s{t} \cdot x^t
& = \frac{1}{1 + V\cdot x + \w{2} \cdot x^2 + \dots + \w{n} \cdot x^n} \\
& = \sum_{p=0}^\infty (-1)^p 
     \bigl(V\cdot x + \w{2} \cdot x^2 + \dots + \w{n} \cdot x^n \bigr)^p.\\
   \end{split}
\end{equation}
For $t \geq 1$, the coefficient of $x^t$ in 
$(V\cdot x + \w{2} \cdot x^2 + \dots + \w{n} \cdot x^n)^p$ is zero
unless $1 \leq p \leq t$, in which case this coefficient is
the class of $T_t^p$.  This proves~\eqref{equation4}.

We now define the differential $\delta$ of~\eqref{equation2} on each
direct summand $\w{i_1} \tensor \cdots \tensor \w{i_p}$ of $T_t^p$, where
$t = i_1 + \dots + i_p$.  We start with the case $p=1$. Given
$v_1, \dots, v_i \in V$, we introduce the notations
$v_A$, for nonempty $A \subset \{1, \dots, i\}$,
and $s(A,B)$, for disjoint nonempty $A, B \subset \{1, \dots, i\}$, by
\begin{equation}
\label{equation8}
\begin{split}
   v_A = v_{a_1} \wedge \dots \wedge v_{a_j} \in \w{j},
     &\text{ where } A = \{a_1, \dots, a_j\}
     \text{ with } a_1 < \dots < a_j,
\\
   s(A,B) \in \{\pm 1\}
     &\text{ such that }
     v_{A \union B} = s(A,B) v_A \wedge v_B.
\end{split}
\end{equation}
Note that $s(A,B)$ depends only on the sets $A$ and~$B$, and not on the
choice of $v_1, \dots, v_i$.  We then define, for a decomposable
tensor $v_1 \wedge \dots \wedge v_i \in \w{i} = T_i^1$,
\begin{equation}
\label{equation9}
 \delta(v_1 \wedge \dots \wedge v_i)
  = \sum_{\substack{\emptyset \neq A,B \subset \{1, \dots, i\} \\
                    A \union B = \{1, \dots, i\}\\
                    A \intersect B = \emptyset\\
          }} s(A,B) v_A \tensor v_B \in T_i^2.
\end{equation}
For example, $\delta(v_1) = 0$, 
$\delta(v_1 \wedge v_2) = v_1 \tensor v_2 - v_2 \tensor v_1$,
and $\delta(v_1 \wedge v_2 \wedge v_3 \wedge v_4)$ is
\begin{equation}
\label{equation10}
\begin{split}
& (v_1 \wedge v_2 \wedge v_3) \tensor v_4
  - (v_1 \wedge v_2 \wedge v_4) \tensor v_3
  + (v_1 \wedge v_3 \wedge v_4) \tensor v_2
  - (v_2 \wedge v_3 \wedge v_4) \tensor v_1\\
& + (v_1 \wedge v_2) \tensor (v_3 \wedge v_4)
  - (v_1 \wedge v_3) \tensor (v_2 \wedge v_4)
  + \text{[3 more terms]}
  + (v_3 \wedge v_4) \tensor (v_1 \wedge v_2)\\
& + v_1 \tensor (v_2 \wedge v_3 \wedge v_4)
  - v_2 \tensor (v_1 \wedge v_3 \wedge v_4)
  + v_3 \tensor (v_1 \wedge v_2 \wedge v_4)
  - v_4 \tensor (v_1 \wedge v_2 \wedge v_3).
\end{split}
\end{equation}
We claim that $\delta$ is well defined, i.e.,
that the right hand side of~\eqref{equation9} is an alternating form in the
vectors $v_1, \dots, v_i$.  The easiest way to verify this is to check that
if two adjacent vectors $v_\ell, v_{\ell+1}$ are equal, then the 
right hand side of~\eqref{equation9} vanishes.

We now define the general action of $\delta$ on 
$\alpha_1 \tensor \cdots \tensor \alpha_p \in 
\w{i_1} \tensor \cdots \tensor \w{i_p} \subset T_t^p$,
for $1 \leq p \leq t-1$, by
\begin{equation}
\label{equation11}
\begin{split}
 \delta(\alpha_1 \tensor \cdots \tensor \alpha_p)
 & = \delta(\alpha_1) \tensor \alpha_2 \tensor \cdots \tensor \alpha_p
  \> - \>  %%%% spacing hack!
     \alpha_1 \tensor \delta(\alpha_2) \tensor \cdots \tensor \alpha_p
\\
  & + \dots + (-1)^{p-1}
    \alpha_1 \tensor \alpha_2 \tensor \cdots \tensor \delta(\alpha_p)
          \in T_t^{p+1}. \\
\end{split}
\end{equation}
Finally, the last differential $\pi: T_t^t \to \s{t}$ is
the natural projection from $T_t^t = V^{\tensor t}$.

\begin{lemma}
\label{lemma1}
The differentials $\delta: T_t^p \to T_t^{p+1}$ and $\pi$ as defined above
are $G$-homomorphisms that satisfy $\delta\delta = 0$ and $\pi\delta = 0$.
In other words, \eqref{equation2} is a complex of $G$-representations.
\end{lemma}
\begin{proof}
From the definition of $\delta(v_1\wedge \dots \wedge v_i)$, we see that it
respects the $G$-action.  (This would have been less transparent if we had
defined $\delta$ from the beginning in terms of basis elements 
$e_{a_1} \wedge \dots \wedge e_{a_i}$ of $\w{i}$, for a fixed basis
$\{e_1, \dots, e_n\}$ of $V$.)  This implies the $G$-linearity in general.
As for $\delta\delta = 0$, one first checks directly that 
$\delta\delta(v_1\wedge \dots \wedge v_i) = 0$.  The crucial ingredient is
that if $A$, $B$, and $C$ are disjoint nonempty subsets of
$\{1, \dots, i\}$, then
$s(A \union B, C) s(A, B) = s(A, B \union C) s(B, C)$.  This
follows from comparing $v_A \wedge v_B \wedge v_C$ with
$v_{A \union B \union C}$.  This settles the case $p=1$, and we then
proceed inductively for larger $p$.  In particular, writing
$\beta = \alpha_2 \tensor \cdots \tensor \alpha_p$, we have
$\delta\delta(\alpha \tensor \beta)
 = \delta\bigl[ \delta(\alpha) \tensor \beta - \alpha \tensor \delta(\beta)
         \bigr]$,
and the reader should be careful to note the $+$~sign in expanding the
first term:
$\delta\bigl[ \delta(\alpha) \tensor \beta ]
 = \delta(\delta(\alpha)) \tensor \beta
 + \delta(\alpha) \tensor \delta(\beta)$,
because $\delta(\alpha)$ has partial degree~$2$.  This proves our result
except at the last step, involving $\pi: T_t^t \to \s{t}$, where one can
show directly that $\pi\delta = 0$.
\end{proof}
\begin{theorem}
\label{theorem2}
The complex~\eqref{equation2} is exact if $t \geq 1$.
\end{theorem}
\begin{proof}
It is more convenient to show instead that the truncated complex
\begin{equation}
\label{equation12}
T_t^{\bullet}: \quad
0 \to T_t^1 \xrightarrow{\delta} T_t^2 \xrightarrow{\delta} \cdots
  \xrightarrow{\delta} T_t^t \to 0
\end{equation}
has zero cohomology everywhere except at the last term, where the
cohomology is $\s{t}$.  The proof is by induction on $t$, the case $t=1$
being trivial.  For the inductive step, we introduce
a filtration on $T_t^{\bullet}$, and compute $H^*(T_t^{\bullet})$ by the
spectral sequence of the filtered complex (see for instance
Proposition~XX.9.1 of~\cite{LangAlgebra}, but note that the indices on the
filtration there run opposite to ours; Section~XI.3 of~\cite{MacLane} is
another reference).  We define our filtration on $T_t^\bullet$ by taking
$\F^f(T_t^p) \subset T_t^p$ to be
\begin{equation}
\label{equation13}
\F^f(T_t^p)
   =  \Directsum_{\substack{i_1 \leq f \\
                            1 \leq i_1, \dots, i_p \leq n\\
                            i_1 + \dots + i_p = t\\}}
             \bigl( \w{i_1} \tensor \cdots \tensor \w{i_p} \bigr).
\end{equation}
Thus $f$ is a bound on the degree $i_1$ of the
``first term'' $\alpha_1$ in any tensor  
$\alpha_1 \tensor \cdots \tensor \alpha_p \in \F^f(T_t^p)$.
We have $T_t^p = \F^t(T_t^p) \supset \dots \supset \F^0(T_t^p) = 0$.
The reader is encouraged at this point to look ahead to the first diagram
in Example~\ref{example3}, which illustrates the complex and
its filtration in the case $t=4$.

The $E_0$ term in our spectral sequence is given by the associated graded
complex $\G(T_t^{\bullet}) = \directsum_{f=1}^t \G^f(T_t^{\bullet})$, and
the spectral sequence abuts to the cohomology $H^*(T_t^{\bullet})$ that we
wish to compute.
Since $\G^f(T_t^p) = \F^f(T_t^p)/\F^{f-1}(T_t^p)$, we have
\begin{equation}
\label{equation14}
  \G^f(T_t^p) 
     \isomorphic 
      \Directsum_{\substack{i_1 = f \\
                            i_1 + \dots + i_p = t\\}}
             \bigl( \w{i_1} \tensor \cdots \tensor \w{i_p} \bigr)
    \isomorphic
        \w{f} \tensor T_{t-f}^{p-1},
\quad\text{for }
1 \leq f \leq t.
\end{equation}
Here we slightly abuse notation, since
$\G^t(T_t^1) = \w{t}$, corresponding to taking $T_0^0 = k$ in the
rightmost term of~\eqref{equation14} when $p=1$ and $f=t$; if $p=1$ and
$f \neq t$, we take $T_{t-f}^0$ and $\G^f(T_t^1)$ to be zero, consistently
with~\eqref{equation3}.
Now $\delta$ respects the filtration, and descends to a differential
$\overline{\delta}: \G^f(T_t^p) \to \G^f(T_t^{p+1})$ that can be
identified with $1 \tensor (-\delta): 
  \w{f} \tensor T_{t-f}^{p-1} \to \w{f} \tensor T_{t-f}^p$.
This is illustrated in the second diagram in Example~\ref{example3}.
Applying our inductive hypothesis, and noting that the presence of
$-\delta$ instead of $\delta$ makes no difference, we see that the
$\overline{\delta}$-cohomology of $\G^f(T_t^\bullet)$ is concentrated
in degree $\bullet = t-f+1$, where it is naturally isomorphic to 
$\w{f} \tensor \s{t-f}$ (this holds even if $f=t$).  Hence the
$E_1$ term of our spectral sequence is
\begin{equation}
\label{equation15}
 0 \to \w{t} \tensor \s{0} \to \w{t-1} \tensor \s{1} 
   \to \dots \to \w{1} \tensor \s{t-1} \to 0,
\end{equation}
with differentials induced from the ``portion'' of the original
differentials $\delta$ on
\begin{equation}
\label{equation16}
 0 \to \w{t} \to \w{t-1} \tensor V \to \w{t-2} \tensor V \tensor V
  \to \dots \to \w{1} \tensor V^{\tensor(t-1)} \to 0. 
\end{equation}
This ``portion'' corresponds to considering in~\eqref{equation9} only the
terms where $B$ has cardinality $\abs{B} = 1$.
Note that \eqref{equation16} is not a complex, as the composition of the
``partial'' $\delta$s is not zero, but the induced maps 
in~\eqref{equation15} do give a complex, due to the presence of the
symmetric powers $\s{i}$ instead of the tensor powers $V^{\tensor i}$.
Comparing the differentials with those in, say, Section~XXI.4
of~\cite{LangAlgebra}, we see that \eqref{equation15} is 
the part of the Koszul complex~\eqref{equation5} in total degree~$t$,
after we drop the final terms $\s{*} \to k \to 0$.  (The fact that the
Koszul complex, as written, starts with $\w{n}$, whereas we have started
with $\w{t}$, is immaterial, since the ``missing'' vector spaces are all
zero, either from $\w{i}$ with $i > n$, or from $\s{j}$ with $j < 0$.)
We have $t \geq 2$, so we can ignore the term $k$ in the Koszul complex,
since that term appears only in total degree~$0$.  Hence
we obtain that the $E_2$ term of our spectral sequence is degenerate,
consisting of a single instance of $\s{t}$ in one corner.  Thus $E_2 =
E_\infty$, and we obtain that the complex~\eqref{equation12} has the desired
cohomology.
\end{proof}

\begin{example}
\label{example3}
We include a diagram of the filtered complex $T_4^{\bullet}$ below,
drawn in such a way that we obtain a second quadrant spectral sequence.
Thus the diagram represents the fourth exact sequence of~\eqref{equation1},
omitting the final term $\s{4}$.
The terms of our complex in a given partial degree $p$ lie on a single
SW-NE diagonal.
Note that this is not a double complex, due to the diagonal arrows.  For
example, the differential
$\delta: \w{4} \to
\bigl(V \tensor \w{3} \bigr)
          \directsum \bigl( \w{2} \tensor \w{2} \bigr)
          \directsum \bigl( \w{3} \tensor V \bigr)
$ 
is represented by the three arrows emanating from $\w{4}$ below, and each
arrow corresponds to the terms in one of the three lines of
formula~\eqref{equation10}.

%\medskip
\begin{center}
\begin{tabular}{c|}
\xymatrix{
      &              &                   & V \tensor \w{3} 
                                                \ar[d]      \\
      &              & \w{2}\tensor\w{2}
                            \ar[r] \ar[d]
                                         & V \tensor \bigl[
                                 (\w{2}\tensor V)\directsum(V\tensor\w{2})
                                                     \bigr]
                                                 \ar[d]    \\
\w{4} \ar@/^4ex/[rrruu]
      \ar[rru]
      \ar[r]
      &  \w{3}\tensor V
        \ar[r]\ar[rru]
                     & \w{2}\tensor V \tensor V
                           \ar[r]        & V\tensor V\tensor V\tensor V \\
}\\
\hline
\end{tabular}
\end{center}
%\medskip
The filtration corresponds to taking all columns to the right of some vertical
line in the above diagram.  Hence each part $\G^f(T_4^\bullet)$ of the
graded complex $E_0$ sees \emph{only} the vertical arrows in a single
column.  Note the 
``common factors'' of $\w{2}$ and of $\w{1} = V$ in the rightmost two
columns above.  We also see trivial ``common factors'' of $\w{4}$ and
$\w{3}$ in the leftmost two columns.  Identifying the remaining factors in
terms of the previous complexes $T_t^{\bullet}$ for $t \leq 3$, we
can represent $E_0$ as the diagram
%\medskip
\begin{center}
\begin{tabular}{c|}
\xymatrix{
      &              &                   & V \tensor T_3^1
                                         \ar[d]_{\overline{\delta}} \\
      &              & \w{2}\tensor T_2^1
                      \ar[d]_{\overline{\delta}}
                                         & V \tensor T_3^2
                                         \ar[d]_{\overline{\delta}} \\
 \w{4}
      &  \w{3}\tensor T_1^1
                      & \w{2}\tensor T_2^2
                                         & V \tensor T_3^3 \\
}\\
\hline
\end{tabular}
\end{center}
%\medskip
We now compute $E_1$.  From the inductive hypothesis, 
only the bottom row survives, with differentials induced from the bottom
row of our first diagram for $T_4^{\bullet}$:
%\medskip
\begin{center}
\begin{tabular}{c|}
\xymatrix{
\w{4} \ar[r]
      &  \w{3}\tensor \s{1} \ar[r]
                      & \w{2}\tensor \s{2} \ar[r]
                                         & V \tensor \s{3} \\
}\\
\hline
\end{tabular}
\end{center}
%\medskip
This is the truncated Koszul complex of~\eqref{equation15}, and hence
computing $E_2$ leaves us with just $\s{4}$ in the SE corner.
\end{example}

\begin{remark}
\label{remark4}
We conclude this note by showing the equivalence
between~\eqref{equation2} and the construction in Sections~2 and~4
of~\cite{Totaro} via the normalized bar resolution.  The normalized bar
resolution is a resolution of $k$ as a left module over an associative
$k$-algebra $\Lambda$ with unit, equipped with an augmentation map 
$\eps: \Lambda \to k$.  We write $\Lbar = \Lambda/k$ for the quotient
of vector spaces, which can if needed be identified with the ideal
$\ker \eps$, although it is preferable not to view $\Lbar$ as a module
over $\Lambda$.  Then the normalized bar resolution of $k$ is
\begin{equation}
\label{equation20}
  \dots \to
\Lambda \tensor \Lbar \tensor \Lbar \tensor \Lbar
   \xrightarrow{\dee_3}
\Lambda \tensor \Lbar \tensor \Lbar
   \xrightarrow{\dee_2}
\Lambda \tensor \Lbar
   \xrightarrow{\dee_1}
\Lambda
   \xrightarrow{\eps}
k
\to 0.
\end{equation}
Here $\tensor = \tensor_k$, as before, and the $\Lambda$-module
structure is via multiplication on the leftmost factor $\Lambda$ in
each tensor product.  The differential $\dee_p$ for $p\geq 1$ is given
by the following formula, where we write $\lbar \in \Lbar$ for the
image of $\lambda \in \Lambda$ (see Section~X.2
of~\cite{MacLane}, with $C = k$, for the proof that this is well
defined and that~\eqref{equation20} is exact):
\begin{equation}
\label{equation21}
\begin{split}
& \dee_p(\lambda \tensor \lbar_1 \tensor \cdots \tensor \lbar_p)
  = \lambda\lambda_1 \tensor \lbar_2 \tensor \cdots \tensor \lbar_p \\
&\qquad
  {} - \lambda \tensor \overline{\lambda_1 \lambda_2}
                \tensor \lbar_3 \tensor \cdots \tensor \lbar_p
  + \cdots 
  + (-1)^{p-1} \lambda \tensor \lbar_1 \tensor \cdots
                \tensor \overline{\lambda_{p-1}\lambda_p} \\
&\qquad
  {} + (-1)^p \eps(\lambda_p) \lambda \tensor \lbar_1 \tensor \cdots
                \tensor \lbar_{p-1}.\\
\end{split}
\end{equation}
We apply the above resolution~\eqref{equation20} in the case when the
algebra $\Lambda$ is the exterior algebra 
$\Lambda = \wedge (V^*)$ on the dual vector space $V^*$.
Thus we can view $\Lbar$ as $\directsum_{i \geq 1} \wedge^i(V^*)
= \directsum_{1 \leq i \leq n} \wedge^i(V^*)$.  Now apply the functor
$\Hom_\Lambda(-,k)$ and note that $\Hom_\Lambda(\Lambda \tensor M, k)$
is naturally isomorphic to the dual $k$-vector space $M^*$.  We obtain
that $\mathrm{Ext}_\Lambda^\bullet(k,k)$ can be computed from the
complex
\begin{equation}
\label{equation22}
 0 \to k 
\xrightarrow{\dee^*_1 = 0} \Lbar^* 
\xrightarrow{\dee^*_2} \Lbar^* \tensor \Lbar^* 
\xrightarrow{\dee^*_3} \cdots.
\end{equation}
We can identify $\Lbar^*$ with $\directsum_{i \geq 1} \w{i}$.  The
complex~\eqref{equation22} is graded by total degree~$t$, which is
compatible with the differentials $\dee^*$.  The components in degree
$t \geq 1$ of~\eqref{equation22} are exactly our $T_t^p$; recall that,
by our conventions, in order for $T_t^p$ to be nonzero it is necessary
to have $1 \leq p \leq t$.  It remains to compare the resulting
differentials $\dee^*$ with the differentials $\delta$, so as to show
that the part of~\eqref{equation22} in total degree~$t$ is essentially
the same as the complex~\eqref{equation12}. 
\end{remark}

\begin{proposition}
\label{proposition5}
In the above context, the restriction to total degree~$t$ of $\dee^*$
in~\eqref{equation22} is the negative of $\delta$, as defined in
\eqref{equation9} and~\eqref{equation11}.
\end{proposition}
\begin{proof}
We shall show only the key step, namely that $\dee^*_2 = -\delta$ for
$\delta$ as in~\eqref{equation9}.  We leave it to the reader to
subsequently verify that the recurrence relating $\dee^*_{p+1}$ to
$\dee^*_2$ is the same as~\eqref{equation11} for $\delta$.  It is
enough to compare $\dee^*_2(\alpha)$ with $\delta(\alpha)$ for
$\alpha = v_1 \wedge \dots \wedge v_i \in \w{i}
   = (\wedge^i(V^*))^* \subset \Lbar^* \subset \Hom_k(\Lambda, k)$.
More precisely, $\alpha$ annihilates $\wedge^{i'}(V^*)$ for 
$i' \neq i$, and acts on
$v^*_1 \wedge \dots \wedge v^*_i \in \wedge^i(V^*)$ by
\begin{equation}
\label{equation23}
\langle \alpha, v^*_1 \wedge \dots \wedge v^*_i \rangle
  = \det \bigl(\langle v_j, v^*_k\rangle \bigr)_{1 \leq j,k \leq i}.
\end{equation}
Unraveling the definitions, we obtain that $\dee^*_2(\alpha) \in
(\Lbar \tensor \Lbar)^*$ acts on $\lbar_1 \tensor \lbar_2$ by
\begin{equation}
\label{equation24}
\bigl \langle \dee^*_2(\alpha), \lbar_1 \tensor \lbar_2 \bigr \rangle
= \eps(\lambda_1) \langle \alpha, \lbar_2 \rangle
  - \langle \alpha, \overline{\lambda_1 \wedge \lambda_2} \rangle
 + \eps(\lambda_2) \langle \alpha, \lbar_1 \rangle.
\end{equation}
The value depends only on the classes $\lbar_1, \lbar_2 \in \Lbar$,
and it is enough to test the action of $\dee^*_2(\alpha)$ in the
situation when $\lambda_1 = v^*_1 \wedge \dots \wedge v^*_{i_1} \in
\wedge^{i_1}(V^*)$ and  
$\lambda_2 = v^*_{i_1+1} \wedge \dots \wedge v^*_{i_1+i_2} \in
\wedge^{i_2}(V^*)$, with $i_1, i_2 \geq 1$.  Hence we have
\begin{equation}
\label{equation25}
%\bigl\langle  \dee^*_2(\alpha), 
%              (v^*_1 \wedge \dots \wedge v^*_{i_1}) \tensor
%                       (v^*_{i_1+1} \wedge \dots \wedge v^*_{i_1+i_2})
%\bigr\rangle
\bigl\langle \dee^*_2(\alpha), \lbar_1 \tensor \lbar_2 \bigr\rangle
= - \langle \alpha, v^*_1 \wedge \dots \wedge v^*_{i_1 + i_2} \rangle.
\end{equation}
The above quantity vanishes unless $i_1 + i_2 = i$, in which case it
is the negative of the determinant of~\eqref{equation23}.  We now wish
to compare this value to 
\begin{equation}
\label{equation26}
\bigl\langle \delta(\alpha), \lbar_1 \tensor \lbar_2 \bigr\rangle
= \sum_{A,B} s(A,B) 
     \langle v_A, \lbar_1 \rangle \cdot \langle v_B, \lbar_2 \rangle
\end{equation}
with $A$ and~$B$ as in~\eqref{equation9}.  This is nonzero only if
there exist choices with $\abs{A} = i_1$ and $\abs{B} = i_2$, which
forces $i = i_1 + i_2$; in that case, the resulting sum
in~\eqref{equation26} reduces to the general Laplace expansion of 
$\det\bigl(\langle v_j, v^*_k\rangle\bigr)_{1 \leq j,k \leq i}$ in
terms of minors of sizes $i_1$ and $i_2$, thereby completing our
proof.
\end{proof}
% include reference to Chapter 7 of Jacobson, BAI?

Our results thus give a direct calculation of the cohomology
of~\eqref{equation22}, and hence of $\mathrm{Ext}^t_\Lambda(k,k)$, in
terms of our calculation of the cohomology of~\eqref{equation12}.
Thus we obtain as a consequence a direct proof of~\eqref{equation0}, using
only the exactness of the Koszul complex.

%%!!!!!!!!!!!!!!!!!!!!!!!!!!!!!!!!!!!!!!!!!!!!!!!!!!!!!!!!!!!!!!!!!!!!!!!!!
%%!!!!!!!!!!!!!!!!!!!!!!!!!!!!!!!!!!!!!!!!!!!!!!!!!!!!!!!!!!!!!!!!!!!!!!!!!
%%!!!!!!!!!!!!!!!!!!!!!!!!!!!!!!!!!!!!!!!!!!!!!!!!!!!!!!!!!!!!!!!!!!!!!!!!!

%%%%%%%%%%%%%%%%%%%%%%%%%%
%\bibliographystyle{amsalpha}
%
%\bibliography{article}
%
%\end{document}
%%%%%%%%%%%%%%%%%%%%%%%%%%%

\providecommand{\bysame}{\leavevmode\hbox to3em{\hrulefill}\thinspace}
\providecommand{\MR}{\relax\ifhmode\unskip\space\fi MR }
\providecommand{\MRhref}[2]{%
  \href{http://www.ams.org/mathscinet-getitem?mr=#1}{#2}
}
\providecommand{\href}[2]{#2}

%%!!!!!!!!!!!!!!!!!!!!!!!!!!!!!!!!!!!!!!!!!!!!!!!!!!!!!!!!!!!!!!!!!!!!!!!!!
%%!!!!!!!!!!!!!!!!!!!!!!!!!!!!!!!!!!!!!!!!!!!!!!!!!!!!!!!!!!!!!!!!!!!!!!!!!
%%!!!!!!!!!!!!!!!!!!!!!!!!!!!!!!!!!!!!!!!!!!!!!!!!!!!!!!!!!!!!!!!!!!!!!!!!!

\end{document}